\documentclass{article}
\usepackage{PRIMEarxiv}

\usepackage[utf8]{inputenc} % allow utf-8 input
\usepackage[T1]{fontenc}    % use 8-bit T1 fonts
\usepackage{hyperref}       % hyperlinks
\usepackage{url}            % simple URL typesetting
\usepackage{booktabs}       % professional-quality tables
\usepackage{amsfonts}       % blackboard math symbols
\usepackage{nicefrac}       % compact symbols for 1/2, etc.
\usepackage{microtype}      % microtypography
\usepackage{lipsum}
\usepackage{fancyhdr}       % header
\usepackage{graphicx}       % graphics
\graphicspath{{media/}}     % organize your images and other figures under media/ folder

\usepackage{tabularx} % We added
\usepackage{multirow} % We added

\usepackage{algorithm} % we added!
\usepackage{algpseudocode}
\usepackage{xcolor}  % we added!

\usepackage{mathtools} % we added!
\usepackage{amsmath,amssymb,amsfonts} % we added!
\usepackage{setspace} % we added

\usepackage{multicol} %we aded

\title{Sensitivity Analysis and Monte Carlo Based Uncertainty Quantification of the In-process Modal Parameters in Milling 

}

\author{
  M. Hashemitaheri \\
  Dept. of Mechanical Engineering \& Engineering Science \\
  University of North Carolina at Charlotte \\
  Charlotte\\
  \texttt{mary.hashemitaheri@gmail.com} \\
   \And
  T.T. Le \\
  Dept. of Mathematics \& Statistics \\
  University of North Carolina at Charlotte \\
  Charlotte \\
  \texttt{tle55@charlotte.edu} \\
  \And
  T. Khan \\
  Dept. of Mathematics \& Statistics \\
  University of North Carolina at Charlotte \\
  Charlotte\\
   \texttt{taufiquar.khan@charlotte.edu}\\
    \And
  H. Cherukuri \\
  Dept. of Mechanical Engineering \& Engineering Science \\
  University of North Carolina at Charlotte \\
  Charlotte\\
  \texttt{harish.cherukuri@charlotte.edu} \\
}

\begin{document}
\maketitle

\begin{center}  \today \vspace{20pt} \end{center}

\begin{abstract}
The material removal rates during milling operations are affected by the selection of the cutting depth and spindle speed. Poor selection of these parameters can result in chatter or suboptimal material removal rates. Stability Lobe Diagrams (SLDs) are the well-known approach to selecting appropriate chatter-free values for these parameters. The Physics-based stability lobe diagram is usually generated using the structural dynamics and the cutting parameters. However, since the machine dynamics are measured in the static state of the machine (zero speed), the generated SLD is not reliable as the machine behavior may vary during the cutting operations. Besides, measuring structural dynamics parameters under cutting conditions is difficult and needs new equipment. This study proposes a new approach to determining in-process structural dynamics parameters based on a multivariate Newton-Raphson method. The physics-based model is combined with empirical records to extract reliable structural dynamics parameters inversely. Some examples based on synthetic data are presented to illustrate this inverse approach. Also, the performance of the algorithm is evaluated on an empirical data set and its ability to improve the stability boundary is verified. Furthermore, the sensitivity analysis is performed to quantify the exposure of the SLD to the changes in each structural dynamics parameter.
\end{abstract}

% keywords can be removed

\setstretch{1.2}

% \begin{multicols}{2}

\section{Introduction}

\label{sec:background}

In manufacturing operations, increasing the material rate removal is an important task to increase productivity. Achieving this goal is possible by increasing spindle speed and depth of cut values. However, the performance of the cutting process can be restricted by the occurrence of unstable vibrations, called chatter. Chatter is harmful as it increases costs in machining operations by damaging the workpiece, tool, and machine. Therefore, many studies have been conducted to avoid chatter and determine the stability properties during cutting operations. Regenerative chatter which is known as the most common type of chatter happening during chip generation is introduced by \cite{tobias1965machine, tlusty1954self}. As a result of involving regenerative effect in the dynamics of milling operation, the equations of motion are delay differential. Many researchers utilized frequency-domain and time-domain methods to analyze the stability properties by generating Stability Lobe Diagram (SLD). The stability lobe diagram is a 2D graph that distinguishes the stable and chatter regions by a stability boundary.
%, as shown in Fig.~\ref{fig:SLD_example}.
The boundary specifies the maximum stable depth of cut for each spindle speed, exceeding which will result in chatter. This diagram enables us to increase material rate removal by choosing proper values of spindle speed and depth of cut while avoiding chatter. The average tooth angle approach is introduced by Tlusty \cite{tlusty1983stability,smith1991overview, smith1990update} as an analytical solution in the frequency domain to obtain the stability boundary. The single frequency solution and the multi-frequency solution are introduced by Altintas and Budak \cite{altintacs1995analytical, budak1998analytical}  to analyze the stability of the milling process using the Fourier series approach. In time-domain solutions, methods like semi-discretization, full-discretization, spectral element, Chebyshev collocation \cite{insperger2011semi,ding2010full, khasawneh2011spectral,totis2014efficient} are utilized to predict chatter with more details included in the simulations. 

Generally, the stability boundary derived based on either of the physics-based models is used to avoid the unstable cut. To generate the SLD, the machine dynamics are identified by impact hammer test at the tool tip in the idle state of the machine. However, in reality, the experimental results do not follow the theoretical boundary built upon the information gathered in the static machine condition. That is, choosing the depth of cut and the spindle speed derived from the static state, may lead to some level of inaccuracy. The discrepancy is due to the changes in machine dynamics during the cutting process. Several studies have been conducted to calculate the in-process structural dynamics parameters to enhance the stability boundary.  Frequency Response Function (FRF) at different spindle speeds by the application of the Operational Modal Analysis (OMA) technique is calculated in \cite{paliwal2020prediction, zaghbani2009estimation}. However, their methods still require an accelerometer to measure the vibration at the tip of the tool. {\"O}z{\c{s}}ahin et al. \cite{ozcsahin2015process} used the inverse stability solution using experimental results to determine the structural dynamics parameters inversely. Their method requires at least two cutting tests. Grossi et al. \cite{grossi2017improved} implemented a spindle speed ramp-up test to quickly extract stability boundary. The authors optimized the in-process structure by minimizing the difference the error between the analytical and experimental chatter frequencies and depth of cut using a Genetic Algorithm. Postel et al.\cite{postel2020neural} built a neural network to identify the cutting coefficient and structural dynamics parameters in milling operations. Their algorithm is trained on experimental cuts to predict underlying parameters. Suzuki et al. \cite{suzuki2012identification} identified the machine behavior by an inverse algorithm that minimizes the error in the depth of cut and phase shifts. Oleaga et al. in \cite{oleaga2022method} proposed a method to improve the frequency response function of milling operations to avoid chatter and increase machining productivity. Many other chatter avoidance experimental techniques and prediction methods are reviewed by Urbikain et al. in \cite{urbikain2019prediction}. 

Many researchers have studied the effect of parameter variations on the stability properties to determine a reliable prediction during machining processes.
Hu et al. \cite{hu2016reliability} applied a neural network to predict the limiting axial depth of cut. The neural network has 9 inputs to generate the stability boundary including cutting force coefficients, modal parameters in each direction (modal damping coefficient, modal stiffness coefficient, modal mass), and radial immersion ratio. The Monte Carlo simulation method and the moment method are utilized to compute the reliability of chatter stability during the milling process. 
Kurdi \cite{kurdi2005numerical} calculated the effect of uncertainty in axial depth of cut on material rate removal and surface location error by evaluating the sensitivity of the maximum stable axial depth of cut at each spindle speed to each of the input parameters (cutting force coefficients, tool modal parameters, and cutting parameters). In their study, time finite element analysis was used to model the milling dynamics. The effect of uncertainties of noisy measured FRF and fitted modal parameters on the stability lobe diagram for frequency-domain and time-domain methods is studied by Hajdu et al. \cite{hajdu2015sensitivity, hajdu2015effect}. Their studies show that the number of modes involved in the calculation of the stability boundary and non-symmetric FRF considerably affect the stability properties. Zhang et al. \cite{zhang2012numerical} developed a numerical robust stability formulation to optimize the upper and lower bound of the stability lobe diagram to address the uncertainties in modal parameters and cutting force coefficients with constraint conditions. In their work, the time finite element approach (TFEA) is used to calculate the stability boundary.

In this study, first, the in-process structural dynamics parameters are extracted, given in-process SLD. Conventionally, the SLD is built having the structural dynamics parameters using either frequency-domain or time-domain methods. However, the calculated SLD is not reliable as these parameters change during cutting operations, and the experimental results do not follow the stability boundary. To overcome this issue, an iterative algorithm is introduced using a multivariate Newton-Raphson method to calculate the in-process structural dynamics parameter. In this approach, the physics-based process model is combined with the experimental results. Through the algorithm, the physics-based SLD and experimental data are compared to each other in each iteration. One may assume that the cutting parameters used to generate the experimental data are fixed. So only the structural dynamics parameters are the input of the algorithm and updated in each iteration. The updated parameters are quantified at the end of each iteration by calculating the error. The algorithm keeps iterating till the updated parameters can explain the experimental results. Three examples based on synthetic data are presented to illustrate the performance of the approach. Also, the algorithm is evaluated on an empirical data set and its ability to improve the stability boundary is verified. Furthermore, results from sensitivity analyses performed to calculate the exposure of the stability lobe diagram to each structural dynamics parameter are presented. Although there are some studies in the literature to calculate in-process structural dynamics parameters, they only focus on a few points on the boundary. This work covers a wide range of spindle speed values by using more data points on the border. Note that in this study, the focus is solely on milling. However, the same method can be applied to other types of machining, e.g., turning. The only difference lies in the machining dynamics model, which varies slightly. Utilizing the right dynamics model, the same inverse method is applicable to any machining operation.

Note that while the presented inverse method provides the predictions for the in-process parameters, the sensitivities quantified on SLD determine to what extent the stability boundary will change in case there are any errors in predictions. The details of the algorithm and sensitivity analysis are explained in Sec.~\ref{sec:2}, followed by the results and discussion in Sec.~\ref{sec:3} and the conclusion in Sec.~\ref{sec:4}.

% \begin{figure}[t]\vspace*{4pt}
% \centering
% \includegraphics[width = 0.65\textwidth]{SLD.png}
% \caption{An example of stability lobe diagram. The boundary distinguishes chatter and stable areas.}
% \label{fig:SLD_example}
% \end{figure}

\begin{table}
    \caption{Values used to generate synthetic data set for the three different examples.}
\label{tab:Input_values}
    \centering
    \resizebox{0.65\textwidth}{!}{%
    \begin{tabular}{lrrrr}
    % \begin{tabularx}{0.48\textwidth}{Xrrr}
    \hline
    Parameter & Ex.1  & Ex.2 & Ex.3    \\ \hline
    Natural frequency: $f_{nx}$   ($Hz$)     & 903.0   & 500.0    & 900.0    \\
    Stiffness: $k_x$               ($N/m$)    & 12.53   & 8.00    & 9.000    \\
    Damping ratio: $\xi_x$                    & 0.0300  & 0.0200  & 0.0200  \\
    Natural frequency: $f_{ny}$    ($Hz$)     & 903.0   & 500.0   & 950.0    \\
    Stiffness: $k_y$               ($N/m$)    & 12.53   & 8.00    & 10.00    \\
    Damping ratio: $\xi_y$                    & 0.0300  & 0.0200  & 0.0100  \\
    Tangential cutting force: $k_t$ ($N/mm^2$)& 556.31  & 695     & 2173 \\
    Normal cutting force:  $k_n$              & 0.404   & 0.404   & 0.268\\
    Start angle: $\phi_s$ ($deg$)             & 0.000   & 0.000   & 126.9\\
    Exit angle: $\phi_e$  ($deg$)             & 180.0   & 180.0   & 180.0\\
    Number of flute: $N_t$                    & 2       &4        & 3\\
    \hline
    \end{tabular}
    % \end{tabularx}
    }
\end{table}

%%%%%%%%%%%%%%%%%%%%%%%%%%%%%%%%%%%%%%%%%%
\section{Methodology}
\label{sec:2}

According to physic-based models, structural dynamics and cutting parameters are determinants of the stability lobe diagram.  In this study, the Fourier series approach (Zero-Order solution) by Altintas \cite{altintas2012manufacturing} is used to generate the stability boundary. The cutting parameter set is comprised of the number of teeth of the tool, the start and exit angle of the tool, and tangential and radial coefficients. Structural dynamics parameters include natural frequency, stiffness, and damping ratio in each direction. These parameters usually are derived experimentally by calculating Frequency Response Function (FRF) which is the result of performing an impact hammer test at the tip of the tool. As mentioned, a goal of this research is to investigate the possibility of extracting the in-process structural dynamics variables, given the empirical stability lobe diagram. Thus, the cutting parameters are fixed as they are pre-defined while collecting experimental data. So only structural dynamics are considered to be determinants of the stability function. It is assumed that $SLD^{exp}$ indicates the experimental SLD which is made during some experimental milling operations considering a set of pre-defined values for cutting parameters. So $SLD^{exp}$ is the reference SLD, and the goal is to return the underlying structural dynamics parameters utilizing a multivariate Newton-Raphson method approach. It is considered that $n$ records are collected during milling, so the experimental SLD data set includes $n$ data points:

\begin{equation}
        \begin{aligned}
SLD^{exp} &= ~\{(\Omega^{exp}_{1},b_1^{exp}), \cdots, (\Omega^{exp}_{n}, b_{n}^{exp})\} 
\end{aligned}
\label{eq:omega}
\end{equation}

where $\Omega$ is spindle speed and $b$ is the maximum depth of cut at which the cut is stable. Note that here SLD is a 2D diagram of the depth of cut vs. spindle speed. So in the resulting SLD from the Fourier series approach \cite{altintas2012manufacturing}, the effect of the feed rate is assumed to be negligible and only the depth of cut and spindle speed are the parameters of interest. Performing hammer test on the tool tip in the idle state of a machine, one may assume that $m$ structural dynamics parameters are derived by calculating the FRF. Hence, according to physic-based models, $SLD^{exp}$ is a function of $m$ structural dynamics parameters, which are the unknowns of the problem and are to be derived through the Newton-Raphson method approach.
\begin{equation}
SLD^{exp} = SLD^{exp}(p_1, p_2,\cdots,p_m) 
\label{eq:sld}
\end{equation}
Newton-Raphson method initiates with an initial guess. So, one may assume some guessed values for the unknown structural dynamics parameters. $(p_1^*,\cdots, p_m^*)$ indicates the initial guess set. An analytical stability lobe diagram is generated based on an initial guess using the physics-based model. $SLD^{*}$ denotes the analytical data set which includes the $n$ data points on the boundary with the same values of spindle speed in the experimental data set.

\begin{equation}
        \begin{aligned}
SLD^{*} &= ~\{(\Omega^{exp}_{1},b_1^{*}), \cdots, (\Omega^{exp}_{n}, b_{n}^{*})\} 
\end{aligned}
\label{eq:omega2}
\end{equation}
Here, $b^{*}$ is the maximum stable depth of cut calculated using the physic-based model at each specific spindle speed. To evaluate the guess, one may need to compare the results obtained from the theory with the reference; That is, comparing the maximum depth of cut values from $SLD^{exp}$ data set (the reference SLD) with the $SLD^{*}$ data set (SLD utilizing physic-based model based on initial guess). Therefore, one may use Mean Log Absolute Error (MLAE) as an objective function for evaluation:
\begin{equation}
MLAE = MLAE(p_1^*,\cdots,p_m^*) =  \frac{1}{n}\Sigma_{i=1}^{n}\log({1+|b_{i}^{exp} - b_{i}^{*}|}) \label{eq:mlae}
\end{equation}
Then, the guessed value for each unknown parameter is updated using Eqn. \ref{eq:newton}: 
\begin{equation}
p_i^* \leftarrow p_i^*- \frac{MLAE \times S_i}{||S||^2}\times\alpha \hspace{15pt}; \hspace{15pt}  i=1,2,\cdots,m.
\label{eq:newton}
\end{equation}
Here, $\alpha$ controls the pace ($\alpha \in (0,1]$), and one may assume $\alpha = 1$ for simplicity of the notation. $S_i$ denotes the sensitivity of the objective function with respect to each parameter which is calculated as follows:
\begin{equation}
        S_i = \frac{\partial MLAE}{\partial p_i^*} = \frac{MLAE(p_1^*,\cdots,p_i^*+\epsilon ,\cdots,p_m^*) - MLAE}{\epsilon}
        \label{eq:partial}
\end{equation}
where, $\epsilon$ is a small value added to a guessed parameter ($p^*_i$) to derive the partial derivative of the objective function with respect to that parameter.

\begin{figure}[t]\vspace*{4pt}
\centering
\includegraphics[width = 0.55\textwidth]{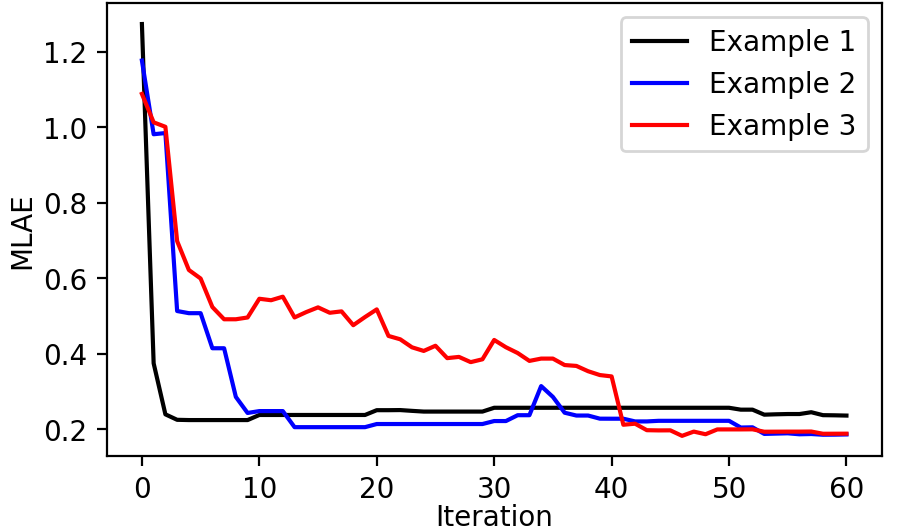}
\caption{The objective function (MLAE) measuring the error terms trough iterations.}
\label{fig:conv}
\end{figure}

\begin{table}
    \caption{Results using Newton method approach for the three examples.}
\label{tab:results_fourier}
    \centering
    \resizebox{0.75\textwidth}{!}{%
    \begin{tabular}{llrrrrrr}
    % \begin{tabularx}{\textwidth}{XXrrrrrrr}
    \hline
 Exp                      & Value        &$f_{nx}$& $k_x$ & $\xi_x$ &$f_{ny}$& $k_y$  &$\xi_y$ \\ \hline
 \multirow{2}{*}{Ex.1}    &  Target      & 903.0  & 12.53 & 0.0300  & 903.0  & 12.53  & 0.0300 \\
                          &  Predicted   & 910.4  & 13.44 & 0.0287  & 910.4  & 13.44  & 0.0286 \\
  \multirow{2}{*}{Ex.2}   &  Target      & 500.0  & 8.00  & 0.0200  & 500.0  &  8.00  & 0.0200 \\
                          &  Predicted   & 503.5  & 8.66  & 0.0196  & 503.5  &  8.66  & 0.0196 \\
 \multirow{2}{*}{Ex.3}    &  Target      & 900.0  & 9.00  & 0.0200  & 950.0  & 10.00  & 0.0100 \\
                          &  Predicted   & 906.8  & 9.49  & 0.0183  & 947.7  &  9.56  & 0.0104 \\ \hline
    \end{tabular}
    }
    % \end{tabularx}
\end{table}

One can repeat the Newton-Raphson step as specified in Eqn. \ref{eq:omega2} to Eqn. \ref{eq:partial} as much as needed. At each iteration, the analytical data set $SLD^{*}$ is generated using the updated values for the parameters. Through this repetitive algorithm, the analytical and experimental stability lobe diagrams are compared to each other at each iteration, and the unknown values are updated at the end of each iteration using the objective function and sensitivity values. Thus, by repeating the algorithm several times, it can converge to the correct answer for each unknown parameter. The algorithm would stop iterating by choosing a proper threshold for the objective function. Hence, the final values extracted through the algorithm could justify the experimental stability lobe diagram.

Also, some adjustments are applied to enhance the Newton method and accelerate its convergence, as dealing with multivariate functions in the Newton method is sometimes challenging. 
In Particular, using the Newton-Raphson method, there is a possibility of getting stuck in the local minima. To overcome this issue, a specific threshold for the objective function is set. If the algorithm was not able to improve the objecting function at least by that threshold after a certain number of consecutive iterations, one may add some e.g. 5\% jumps in each/some parameters. However, adding random jumps does not necessarily improve the results. Thus, the jumps are tried, but only applied if they can successfully improve the objective function. This way the algorithm gets a chance to run away from local minima. Moreover, during the experiments, it is observed that sometimes the algorithm fails on the critical points of the stability lobe diagram. Therefore, the algorithm needs to be more sensitive to the SLD errors at those critical points. That is achievable by giving more weights to the critical points rather than equal weights that are implied in Eqn.\ref{eq:mlae}. Finally, the convergence of the algorithm may depend on the choice of the initial guess. While the algorithm may converge fast using one initial guess, using another initial guess it may need many more iterations to converge. Therefore, a diverse set of initial guesses are utilized in practice. The algorithm initially may try different guesses, but after a few iterations, only continues with the best guess.

\begin{figure*}[t]\vspace*{4pt}
\centering
\includegraphics[width = 1\textwidth]{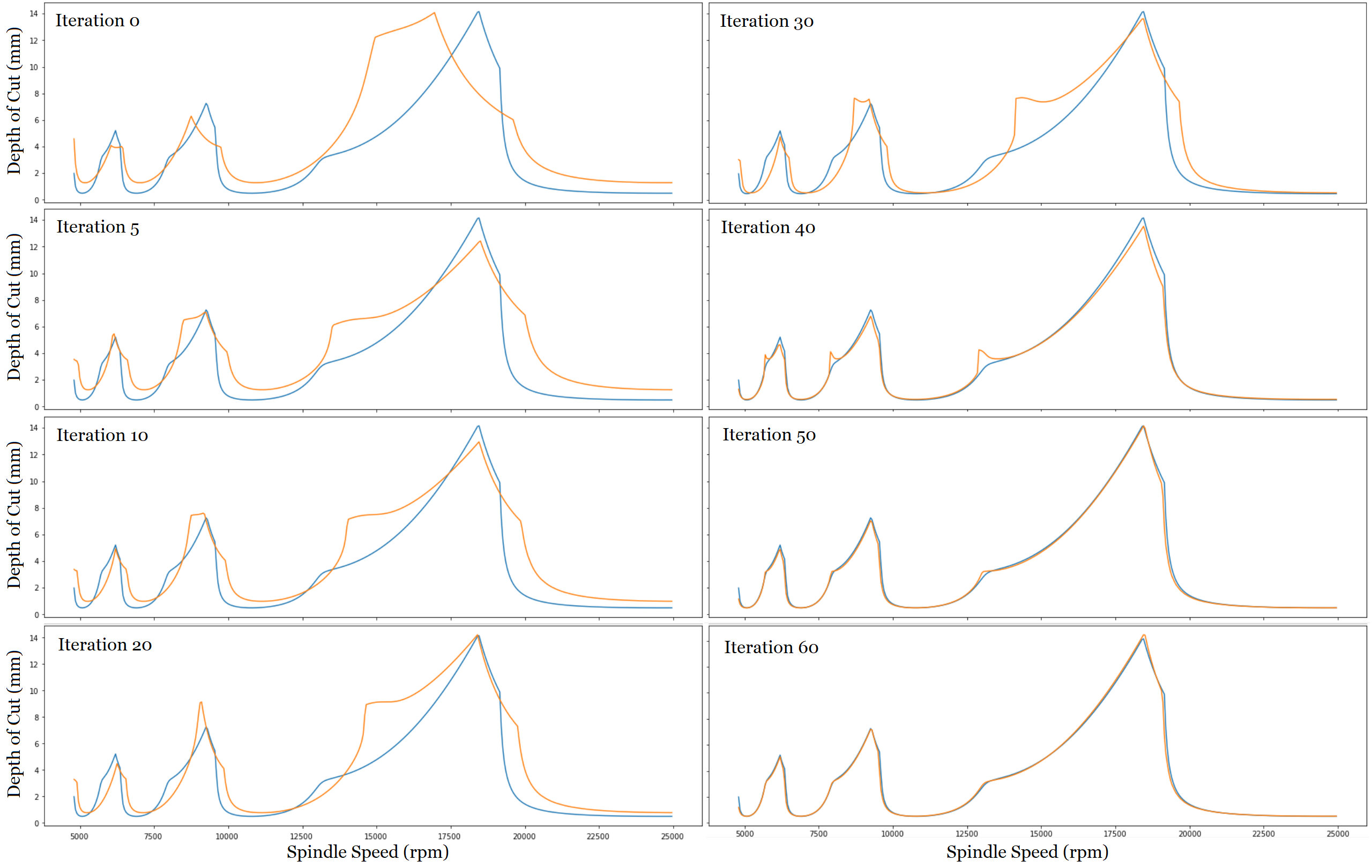}
\caption{Some of the steps in Newton-Raphson approach for example 3. The blue curve is the empirical SLD while the orange curve is the guessed SLD at each step.}
\label{fig:iter}
\end{figure*}

\begin{figure*}[t]\vspace*{4pt}
\centering
\includegraphics[width = 1\textwidth]{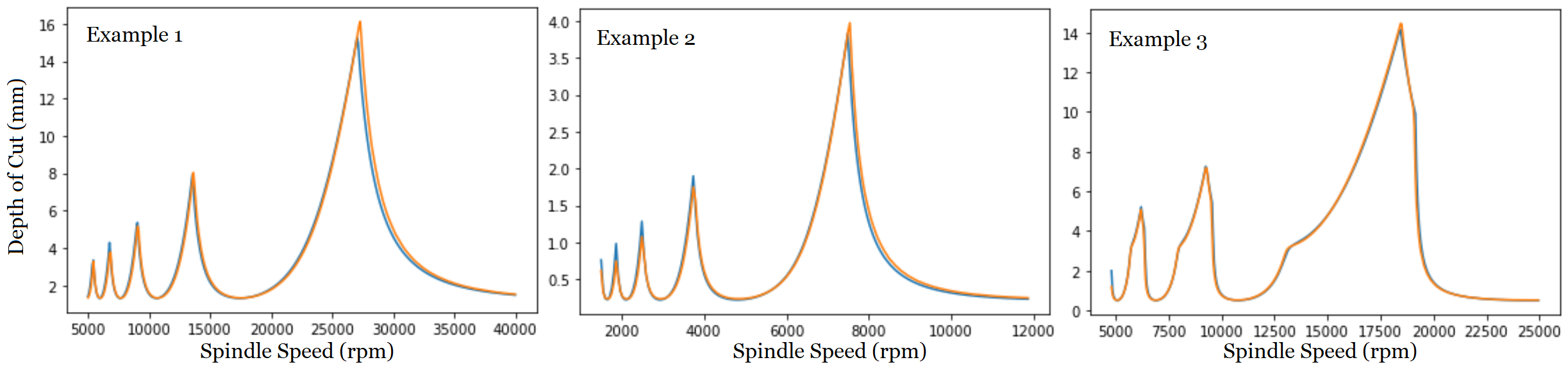}
\caption{Final SLDs for the three examples. Blue curves are the target SLDs while the orange curve shows the calculated SLD in each example.}
\label{fig:eggs}
\end{figure*}

%%%%%%%%%%%%%%%%%%%%%%%%%%%%%%%%%%%%%%%%%%
\section{Results}
\label{sec:3}
This section studies three different examples of SLDs. For each case, first, the Newton-Raphson approach is applied to extract the underlying structural dynamics parameters as explained in section \ref{sec:2}. Then, the results of the sensitivity analysis on SLDs are presented, followed by the discussion on the exposure of the stability lobe diagram to changes in each of the structural dynamics parameters in section \ref{subsec:Sensitivity}.

\subsection{Newton-Raphson Method}
\label{subsec:NM}

Three different examples are studied to derive the structural dynamics parameters using the multivariate Newton-Raphson method approach. For each example, synthetic data is generated using the Fourier series approach (Zero-Order solution \cite{altintas2012manufacturing}). Each data set includes values of spindle speed and the corresponding depth of cut on the stability boundary.
Table \ref{tab:Input_values} shows the structural dynamics parameters and cutting parameters used to generate these data sets.
In each example, the SLD made using the information provided in Table \ref{tab:Input_values} is considered as the reference SLD. The guessed SLD is made through the Newton-Raphson algorithm and then compared to the reference SLD.

It is assumed that the synthetic data set for each case comes from an experiment, and the values for cutting parameters are given, i.e., normal and tangential cutting forces, number of flutes, and start and exit angles. Thus, only the structural dynamics parameters are the unknowns, and the goal here is to derive them through the Newton-Raphson approach. It is assumed that there is only one dominant mode in each direction for each example. Therefore, there exist totally six unknown parameters for each example. The Newton-Raphson approach is applied for each example separately. The SLD generated through the algorithm to be compared with the reference SLD is also built based on the Fourier series approach (Zero-Order solution). The objective function used for the evaluation of parameters in each iteration is mean log absolute error. In each example, the Newton-Raphson method starts with an initial guess for the parameter set. All the modifications and improvements explained in section \ref{sec:2} are applied to the algorithm. The algorithm keeps iterating until it discovers the parameters that can justify the reference SLDs. 

Table \ref{tab:results_fourier} presents the results derived from Newton-Raphson algorithm. The predicted values are the results found using the algorithm and the target values are the parameters used to build the reference SLD. The results show that the algorithm converges and finds the parameters that can justify the reference SLD for all the cases, as the predicted values are fairly close to their targets. Figure \ref{fig:eggs} shows the reference and the predicted stability lobe diagrams. The blue curves are generated based on target values and are called the reference SLDs, and the orange curves are generated based on the predicted values after 60 iterations and are called the guessed SLDs. As shown in Figure \ref{fig:eggs}, in each example, the two boundaries are very close to each other after 60 iterations. Some of the Newton-Raphson method iterations are visualized for example 3 in Figure \ref{fig:iter}. As Figure Figure \ref{fig:conv} and \ref{fig:iter} show, the iterative process does not necessarily imply a monotone convergence. Nevertheless, after enough iterations, the guessed SLD is similar to the reference SLD. Note that one may keep iterating the algorithm more to reduce the errors. However, that may also imply more computational cost, while the errors are already negligible, and the results might be already acceptable.

\begin{figure}[t]\vspace*{4pt}
\centering
\includegraphics[width = 0.5\textwidth]{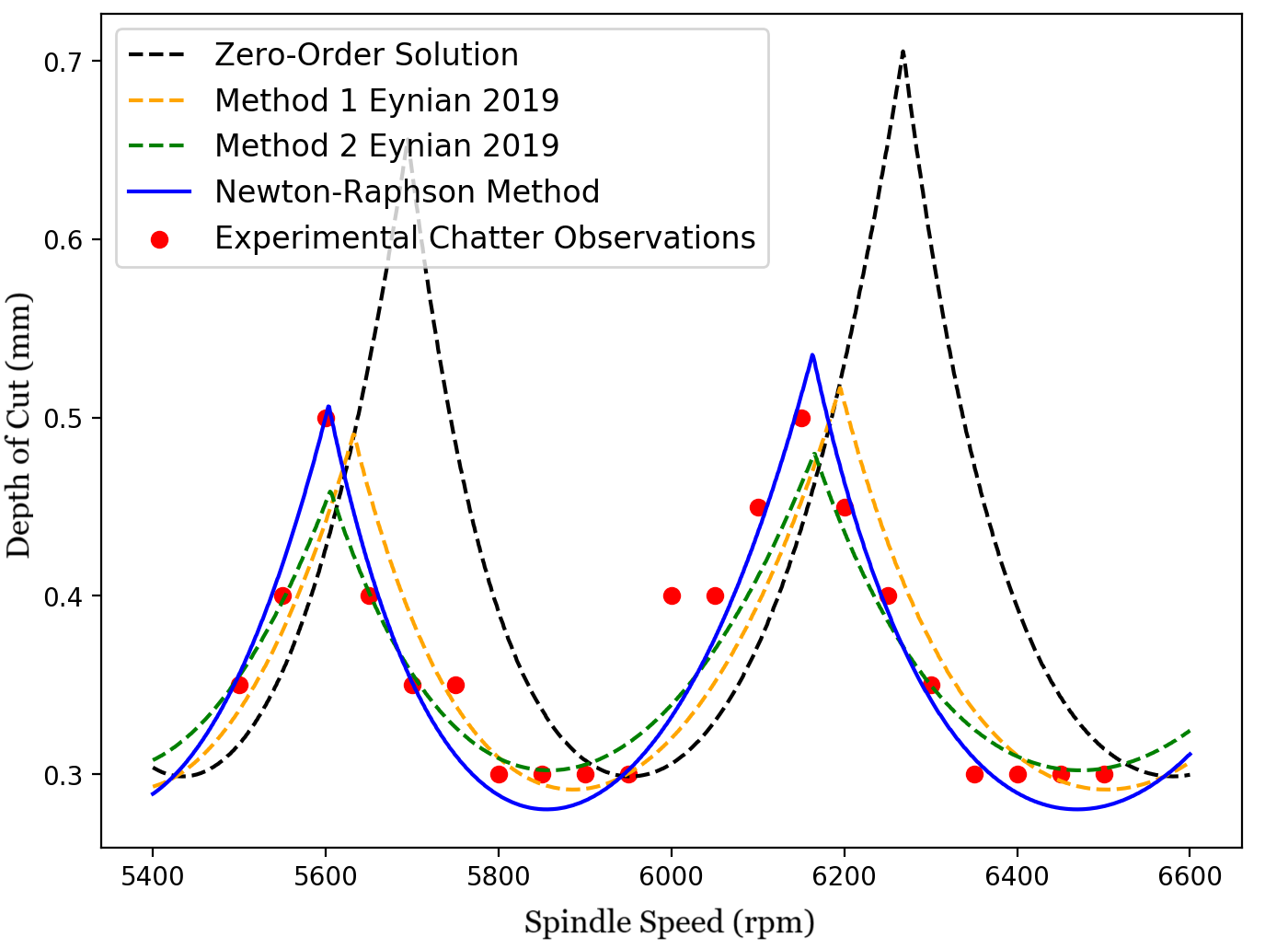}
\caption{Stability lobe diagram on empirical data for Newton-Raphson method and methods used in \cite{eynian2019process}.}
\label{fig:fourSLD}
\end{figure}

To evaluate the Newton method, also an experimental data set presented in \cite{eynianDATA, eynian2019process} is utilized. This experiment was performed with a four-flute tungsten carbide end-mill. The workpiece is a plate made of aluminum 7075-T651. The data set contains $21$ chatter tests of slot milling on the boundary at spindle speed between $5500$ and $6500 rev/min$, which are illustrated by red points in the Figure \ref{fig:fourSLD}. That is, any point below these points is recognized as stable and any point above them  is recognized as chatter. The tangential cutting force, calculated in static conditions is $1110$MPa, and the radial cutting force coefficient is 0.22. The modal parameters measured by applying the impact test in the static state of machine are reported in Table \ref{tab:modal}, where $x$ and $y$ are feed and cross-feed directions, respectively. Stability lobes in zero spindle speed is generated using the information in Table \ref{tab:modal}. This boundary is shown by black dashed curve in Figure \ref{fig:fourSLD}, which predicts many of chatter cuts as stable incorrectly. To capture the in-process structural dynamics parameters, the Newton-Raphson method is applied to the experimental data set. These experimental records are considered as the target SLD, however, there are fewer data points making the boundary comparing the target SLD in synthetic examples. Table \ref{tab:experimentalResults} shows the results derived from Newton-Raphson method. In this case, it is assumed that the tool-spindle assembly is axisymmetric, so the structural dynamics parameters have equal properties in the feed and cross-feed directions. The stability boundary generated using these values are shown in Figure \ref{fig:fourSLD} with blue solid curve. The natural frequency value decreases notably to $4103.4$Hz comparing to the results of the impact test ($4182$ Hz). Thus, the experimental chatter records are shifted to the left-hand side of the stability lobe diagram on stationary speed. The stiffness value reaches to $11.65$ MN/m which is considerably less than the results of the impact test $(15.40MN/m)$ while the damping ratio increases from $0.0170$ in stationary condition to $0.0269$. Table \ref{tab:experimentalResults} also includes the results from two method used in \cite{eynian2019process}(Method 1: Two-point Method, Method 2: Regression Method). As the table suggests, the results are almost consistent with the results extracted from these two methods, but as Figure \ref{fig:fourSLD} shows, the SLD generated based on the results of Newton-Raphson method covers more chatter records.

\begin{table}
    \caption{Structural dynamics parameters in zero spindle speed for the four-flute tool as reported in \cite{eynian2019process}.}
\label{tab:modal}
    \centering
    \resizebox{0.5\textwidth}{!}{%
    \begin{tabular}{lrrrr}
    % \begin{tabularx}{0.48\textwidth}{Xrrr}
    \hline
    Parameter & Mode.1  & Mode.2     \\ \hline
    Natural frequency: $f_{nx}$   ($Hz$)      & 3890.0   & 4182.0       \\
    Stiffness: $k_x$               ($MN/m$)    & 22.60   & 15.40       \\
    Damping ratio: $\xi_x$                    & 0.0196  & 0.0170    \\
    Natural frequency: $f_{ny}$    ($Hz$)     & 3872.0   & 4127.0       \\
    Stiffness: $k_y$               ($MN/m$)    & 23.40   & 25.10        \\
    Damping ratio: $\xi_y$                    & 0.0220  & 0.0177    \\
    \hline
    \end{tabular}
    }
    % \end{tabularx}
 \end{table}

\begin{table}
    \caption{In-process structural dynamics parameters calculated using Newton-Raphson method and the methods used in \cite{eynian2019process}.}
\label{tab:experimentalResults}
    \centering
    \resizebox{0.75\textwidth}{!}{%
    \begin{tabular}{lrrrr}
    % \begin{tabularx}{0.48\textwidth}{Xrrr}
    \hline
    Parameter & Newton-Raphson  & Method.1 & Method.2     \\ \hline
    Natural frequency: $f_{n}$   ($Hz$)      & 4103.4   & 4124   &   4103   \\
    Stiffness: $k$              ($MN/m$)    & 11.65     & 10.63    &  8.90  \\
    Damping ratio: $\xi$                    & 0.0269      & 0.0307  & 0.0379\\
    % Natural frequency: $f_{ny}$    ($Hz$)     & 3872.0   & 4127.0       \\
    % Stiffness: $k_y$               ($MN/m$)    & 23.40   & 25.10        \\
    % Damping ratio: $\xi_y$                    & 0.0220  & 0.0177    \\
    \hline
    \end{tabular}
    % \end{tabularx}
    }
\end{table}

\renewcommand{\algorithmicrequire}{\textbf{Input:}}
\renewcommand{\algorithmicensure}{\textbf{Output:}}

\begin{algorithm}
\caption{Sensitivity analysis on exact parameter set}\label{alg:1}
\begin{algorithmic}
\Require ~~Parameters=$\{p_1, \cdots, p_m\}$
\Ensure MSE in SLD when changing each parameter
\State Make the original SLD: $SLD^*=\{(\Omega_i^*,b_i^*)| ~i=1,\cdots,n\}$

\For {$j$ in $1,\cdots,m$}
    \For{$\epsilon$ in $\{-0.20, -0.15, \cdots, +0.20 \}$}
        \State Change the parameters to $P'$:
        \State ~~~~ $P'=\{p_1, \cdots, (1+\epsilon)p_j, \cdots, p_m\}$
        \State Make modified SLD using $P'$: 
        \State ~~~~ $SLD'=\{(\Omega_i^*,b_i^{'})| ~i=1,\cdots,n\}$
        \State Compute MSE: $MSE_{j,\epsilon} = \frac{1}{n}\sum_{i=1}^{n} (b_i^* - b_i^{'})^2$
    \EndFor
\EndFor
\end{algorithmic}
\end{algorithm}

\begin{figure*}[t]\vspace*{4pt}
\centering
\includegraphics[width = 1\textwidth]{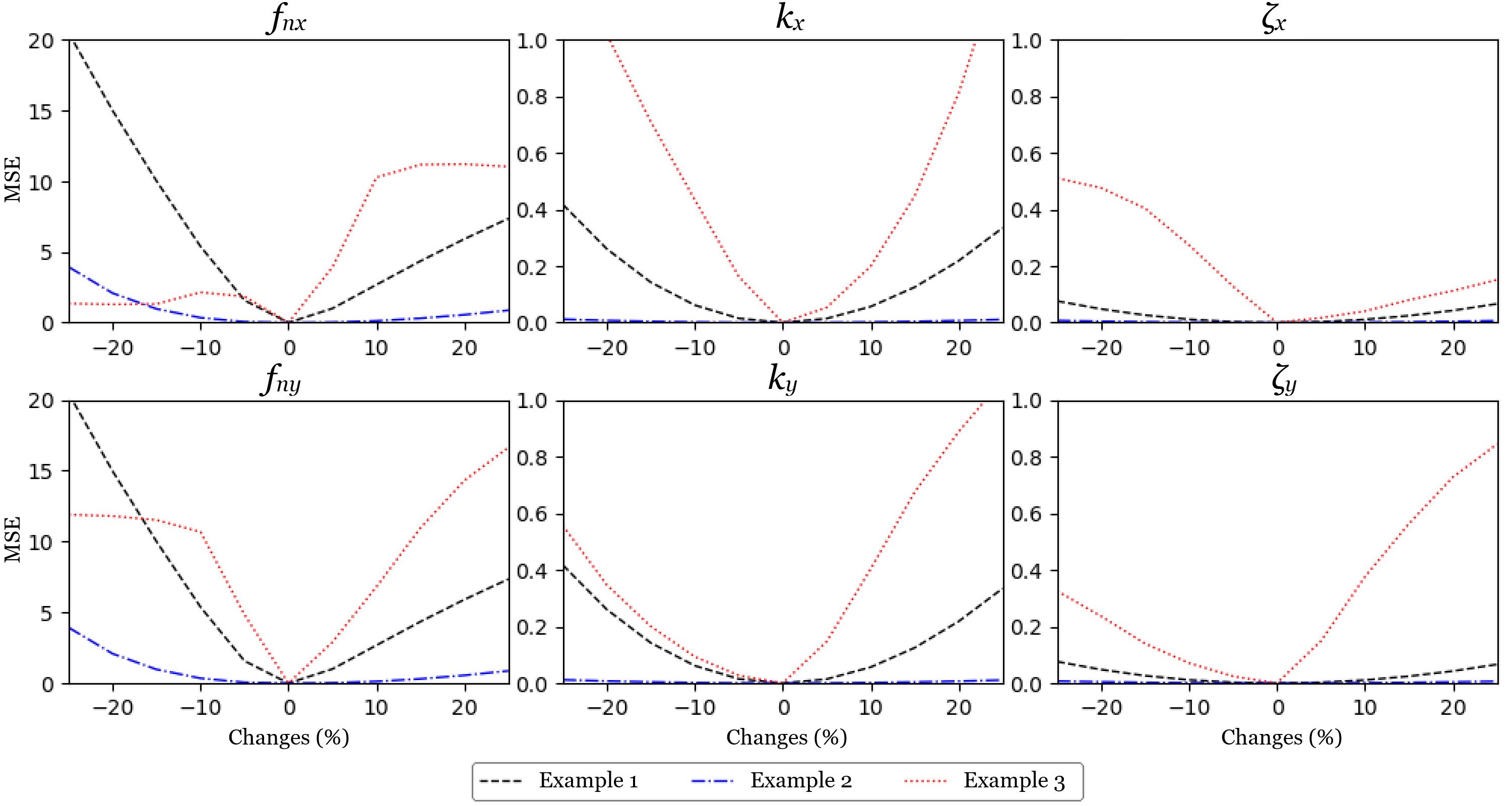}
\caption{Sensitivity of the stability lobe diagram to the changes in structural dynamics parameters.}
\label{fig:SENS}
\end{figure*}

\subsection{Sensitivity Analysis}
\label{subsec:Sensitivity}

To analyze the sensitivity of the stability lobe diagram with respect to different parameters, two different approaches are utilized. First, three different examples of SLD using three different parameter sets are tried. For each SLD, the parameters are changed one-by-one by $\epsilon$\% and a modified SLD is generated for different values of $\epsilon \in [-25, 25]$. Then the modified SLD is compared to the Original SLD using Mean Squared Error (MSE) for different values of $\epsilon$. 
Note that when each parameter is changed, the other parameters are kept constant. The different MSE plots will help us to analyze the sensitivity of the SLD to each parameter, when that parameter is changed by different ratios. The above-mentioned steps are summarized in Algorithm \ref{alg:1}. The predicted parameter sets $(f_{nx}, k_x, \zeta_x, f_{ny}, k_y, \zeta_y)$ as specified in Table \ref{tab:results_fourier} are used for the three examples.

Figure \ref{fig:SENS} visualizes the results. As shown in the figure, the sensitivity of the SLD with respect to different parameters differs from one example to another. However, generally, the SLD appears to be much more sensitive to $f_{nx}$ and $f_{ny}$ than to the other structural dynamics parameters. Note that the vertical scale for these two parameters is 20 times larger than other subplots in Figure \ref{fig:SENS} to cover a different range of changes. 

The variations in sensitivity plots of different examples as reflected in Figure \ref{fig:SENS}, inspire us to use a more generalized approach to aggregate the observations of the sensitivity in different cases with a large number of simulations.

\begin{algorithm}
\caption{Monte Carlo based sensitivity analysis}\label{alg:2}
\begin{algorithmic}
\Require ~~Parameters=$\{p_1, \cdots, p_m\}$, ~Neighborhood width $t$, Number of paths N
\Ensure MSE in SLD when changing each parameter
\For {$path$ in $1,\cdots,N$}:
    \For {$j$ in $1,\cdots,m$}
        \State $upperbound \gets (1+t) \times p_j$
        \State $lowerbound \gets (1-t) \times p_j$
        \State $p_j^* \gets uniform[lowerbound~,~upperbound]$
    \EndFor
    \State Make the reference SLD using $P^*$: 
    \State ~~~~   $SLD^*=\{(\Omega_i^*,b_i^*)| ~i=1,\cdots,n\}$
    
    \For {$j$ in $1,\cdots,m$}
        \State Modify the parameters to $P'$ changing $p_j^*$:
        \State ~~~~ $P'=\{p_1^*, \cdots, (1+\epsilon)p_j^*, \cdots, p_m*\}$
        \State Make modified SLD using $P'$: 
        \State ~~~~ $SLD'=\{(\Omega_i^*,b_i^{'})| ~i=1,\cdots,n\}$
        \State Compute MSE: 
        \State ~~~~ $MSE_{j, path} = \frac{1}{n}\sum_{i=1}^{n} (b^* - b')^2$
    \EndFor
\EndFor

\For {$j$ in $1,\cdots,m$}
    \State Aggregate the results over parameter j:
    \State ~~ $meanMSE_j = \frac{1}{N}\sum_{path=i}^{N} MSE_{j, path}$
    \State ~~ $stdMSE_j ~~~~= \{\frac{1}{N}\sum_{path=i}^{N} (MSE_{j, path} - meanMSEj)^2\}^{1/2}$
\EndFor

\end{algorithmic}
\end{algorithm}

In the second approach, the study is not confined to the specific parameter set $P=(p_1,\cdots,p_m)$ based on which the original SLD is generated. Instead, one may search and measure the sensitivity of SLD in the neighborhood around given parameter set. In a Monte Carlo simulation, one may repeatedly choose a new parameter set $P^*=(p_1^*,\cdots,p_m^*)$  where each parameter $p_i^*$ is randomly chosen from the uniform distribution e.g. from $[0.9p_i ~,~ 1.1p_i ]$ or generally from $[(1-t)p_i ~,~ (1+t)p_i]$ where t defines the width of the neighborhood. Then the sensitivity of the SLD is measured with respect to each corresponding parameter $p_i$, using MSE between the SLD using $P^*$ and its modified version where the corresponding parameter is again changed by a small ratio, e.g. $\epsilon=1\%$. Finally, different results from different Monte-Carlo paths are aggregated over each parameter using aggregators such as mean and standard deviation. The steps are summarized in Algorithm \ref{alg:2}. Table \ref{tab:MC} shows the results of the Monte Carlo simulation with $N=10,000$ iterations aggregating the sensitivities by the means and the standard deviations. The predicted values in example 1 are used as the given parameter set around which the neighborhood of parameters is constructed.

It is immediately observable from Table \ref{tab:MC} that SLD sensitivities show large variations in different simulations. But generally, the SLDs are much more sensitive to $f_{nx}$ and $f_{ny}$ than the other structural dynamics parameters.

\begin{table}[!ht]
    \caption{Results of Monte Carlo simulation measuring the sensitivity of SLD with respect to different structural dynamics parameters.}
\label{tab:MC}
    \centering
    \resizebox{0.85\textwidth}{!}{%
    \begin{tabular}{lcccccc}
    \hline
    Measure & $f_{nx}$ & $k_x$ & $\zeta_x$& $f_{ny}$ & $k_y$ & $\zeta_y$ \\ \hline
    
    Mean Sensitivity  & 0.06610 & 0.00080 & 0.00025 & 0.06798 & 0.00078 & 0.00026 \\
    SD of Sensitivity & 0.03202 & 0.00029 & 0.00012 & 0.03174 & 0.00026 & 0.00013 \\ \hline
    \end{tabular}
    }
\end{table}

%%%%%%%%%%%%%%%%%%%%%%%%%%%%%%%%%%%%%%%%%%
\section{Conclusions}
\label{sec:4}

In this paper, the sensitivity of stability boundary with respect to machine dynamics parameters was investigated. Generally, these parameters are calculated by applying an impact hammer test in the idle state of the machine which results in generating an unreliable stability boundary. This happens because these parameters vary during the cutting process. Thus, a novel approach was introduced to find the in-process structural dynamics parameters under which the physics-based models can support the empirical data. In this approach, a multivariate Newton-Raphson method was utilized that combines the theory with the experimental results. Three synthetic data sets were used to replace the experimental boundary in this study. The algorithm was evaluated as it converged to the correct answers in those three cases and successfully found the underlying structural dynamics parameters. Also, the algorithm was evaluated on an empirical data set and verified its ability to improve the stability boundary. The results derived from the algorithm were used to discover the sensitivity of stability boundary with respect to each parameter. Two different methods were utilized for this purpose. Although the stability border shifts as any structural dynamics parameter changes, the results showed that the natural frequency is the most sensitive parameter. 
In the future, the authors will use more empirical data instead of synthetic data to derive the in-process machine behavior. The synthetic data sets usually suffer from simplifications of the physics-based models. Also, the authors would use faster root-finding methods that can possibly replace Newton-Raphson. While it is not guaranteed that such methods will produce better results for this particular study, using them in the future will not change the proposed problem-specific adjustment in the presented framework.

\section*{Author Contributions}Conceptualization, Maryam Hashemitaheri, Harish Cherukuri and Taufiquar Khan; Formal analysis, Maryam Hashemitaheri; Funding acquisition, Harish Cherukuri; Methodology, Maryam Hashemitaheri; Software, Maryam Hashemitaheri; Supervision, Harish Cherukuri; Writing – original draft, Maryam Hashemitaheri; Writing – review \& editing, Maryam Hashemitaheri, Thuy Le, Harish Cherukuri and Taufiquar Khan.

\section*{Acknowledgments}The Center for Self-Aware Manufacturing and Metrology, and the research described in this paper, are supported under a multi-year grant from the University of North Carolina’s Research Opportunities
Initiative.

\section*{Conflicts of Interest}The authors declare no conflict of interest.

\section*{Note}
Some parts of this article were previously presented orally as \cite{ht2022} and \cite{hashemitaheri2022extracting}.

% \end{multicols}

%Bibliography
\bibliographystyle{unsrt}  
% \bibliography{references} 

\end{document}